\documentclass[a4paper]{amsart}

\usepackage[utf8]{inputenc}
\usepackage[english]{babel}
\usepackage{slashbox,color,enumerate}
\usepackage[pdftex]{hyperref}
\DeclareMathOperator{\Id}{Id}
\newtheorem{theo}{Theorem}

\newtheorem{lemma}[theo]{Lemma}
\newtheorem{prop}[theo]{Proposition}
\newtheorem{cor}[theo]{Corollary}

\newtheorem{assum}[theo]{Assumption}
\newtheorem{assums}[theo]{Assumptions}

\newtheorem{rem}[theo]{Remark}

\newtheorem{exa}[theo]{Example}

\numberwithin{theo}{section}
\numberwithin{equation}{section}

 \def\mG{\mathsf{G}}
 \def\mV{\mathsf{V}}
 \def\mE{\mathsf{E}}

 \def\mv{\mathsf{v}}
 \def\me{\mathsf{e}}

 \def\mf{\mathsf{f}}

\DeclareMathOperator{\DN}{D\!N}

\author{Delio Mugnolo}
\address{Delio Mugnolo, Institut f\"ur Analysis, Universit\"at Ulm, Helmholtzstra{\ss}e 18, 89081 Ulm, Germany}
\email{delio.mugnolo@uni-ulm.de}

\title[Asymptotics of semigroups generated by operator matrices]{Asymptotics of semigroups generated by operator matrices}
\thanks{The author is supported by the Land Baden--Württemberg in the framework of the \emph{Juniorprofessorenprogramm} -- research project on ``Symmetry methods in quantum graphs''.}
\begin{document}

\cleardoublepage
\maketitle

\begin{abstract}
We survey some known results about operator semigroup generated by operator matrices with diagonal or coupled domain. These abstract results are applied to the characterization of well-/ill-posedness for a class of evolution equations with dynamic boundary conditions on domains or metric graphs. In particular, our results on the heat equation with general Wentzell-type boundary conditions complement those previously obtained by, among others, Bandle--von Below--Reichel and Vitillaro--Vázquez.
\end{abstract}

\section{Introduction}

Elliptic problems for second order differential operators with eigenvalue dependent boundary conditions like
\begin{equation}\label{EBP}
\left\{
\begin{array}{rcll}
\Delta u &=& \lambda u\qquad &\hbox{in }\Omega,\\
-\frac{\partial u}{\partial \nu} &=& \lambda u &\hbox{on }\partial \Omega,
\end{array}
\right.
\end{equation}
 have a long history. They appear in certain models of mathematical physics, quantum chemistry, and neuroscience -- cf.\ e.g.~\cite{Esc94} for a list of older references -- but they are probably best known in view of their interplays with certain classes of stochastic processes: Indeed, taking the trace on $\partial \Omega$ of the first equation and plugging it into the second equation we find
\begin{equation}\label{eq:gwbcnong}
\Delta u  = -\frac{\partial u}{\partial \nu}\qquad \hbox{on }\partial \Omega,
\end{equation}
a boundary condition that was first extensively studied by A.\ Wentzell in~\cite{Ven60}.

Another motivation for the study of this boundary condition comes from the observation that the Rayleigh quotient for this problem is
\[
\frac{\int_\Omega |\nabla u|^2 dx}{\int_\Omega | u|^2 dx+\int_{\partial \Omega} | u|^2 d\sigma},
\]
as a direct application of the Gauß--Green formulae shows (here we denote by $\sigma$ the surface measure of $\partial \Omega$). This can be seen as an overlapping of the Rayleigh quotients
\[
\frac{\int_\Omega |\nabla u|^2 dx}{\int_\Omega | u|^2 dx}
\]
and
\[
\frac{\int_\Omega |\nabla u|^2 dx}{\int_{\partial \Omega} | u|^2 d\sigma}\ ,
\]
which by Courant's minimax principle yield the eigenvalues of the Neumann and Steklov problems, respectively.

After studying this elliptic system, it is natural to discuss the Cauchy problem for the corresponding parabolic problem
\[
\left\{
\begin{array}{rcll}
\frac{\partial u}{\partial t}(t,x)&=&\Delta u(t,x) \qquad &t\ge 0,\ x\in \Omega,\\
\frac{\partial u}{\partial t}(t,z)&=&-\frac{\partial u}{\partial \nu}(t,z) \qquad &t\ge 0,\ z\in \partial \Omega,\\
u(0,x)&=&u_0(x),&x\in \Omega,\\
u(0,z)&=&u_1(z),&z\in \partial \Omega,
\end{array}
\right.
\]
say for initial data $u_0,u_1\in L^2(\Omega)\times L^2(\partial \Omega)$. (Again,~\eqref{eq:gwbcnong} is obtained taking the trace of the first equation and plugging it into the second one.)

The fact that (a natural realization of) the Laplacian with boundary condition as in~\eqref{eq:gwbcnong} generates a submarkovian semigroup on $L^2(\Omega)\times L^2(\partial \Omega)$ was first proved in~\cite{AmaEsc96} by Dirichlet form methods and, independently, in~\cite{FavGolGol02} by an application of the Lumer--Phillips Theorem; but the relevant Dirichlet form had already been determined in~\cite{Fuk69}. Actually, there is a long series of different but complementary results on Wentzell-type boundary conditions in the theory of Dirichlet forms, see e.g.~\cite{Pos13} for an updated series of references.

Both papers \cite{AmaEsc96,FavGolGol02} observe that this Hilbert space semigroup extends to a family of \emph{analytic} semigroups on all spaces  $L^p(\Omega)\times L^p(\partial \Omega)$ for $1<p<\infty$. The methods of~\cite{AmaEsc96} were refined in~\cite{AreMetPal03} to prove generation (but no analyticity) in the spaces $L^1(\Omega)\times L^1(\partial\Omega)$ and $C(\overline{\Omega})$: All these papers exploit more or less directly the variational structure of the problem, and in particular the fact that the relevant semigroup is contractive with respect to the norm of $L^\infty(\Omega)\times L^\infty(\partial \Omega)$. However, it is known that even simple lower order perturbations of the Laplacian generate semigroups that are \emph{not} $L^\infty(\Omega)\times L^\infty(\partial \Omega)$-contractive, cf.~\cite[Thm.~4.3]{Nit12}, so that the applicability of interpolation methods becomes less apparent unless one applies more sophisticated methods like in~\cite{AlbHoeStr77}.

Furthermore, additional efforts were needed to prove analyticity of the semigroup on $C(\overline{\Omega})$ and $L^1(\Omega)\times L^1(\partial \Omega)$: This issue could only be settled in \cite{EngFra05,War10} when quite different techniques were first exploited. In particular, in~\cite{EngFra05} the authors applied in a smart way a theory of a certain class of operator matrices that has been developed by K.-J.\ Engel beginning with~\cite{Eng97}.

As we already mentioned, a motivation for the study of eigenvalue dependent boundary conditions comes from the theory of stochastic processes. Indeed, a characterization of all boundary conditions associated with Markov processes is known, see e.g.~\cite{Tai04}, and suggests that additional terms may conveniently be added to~\eqref{eq:gwbcnong} to achieve maximal generality: Following Wentzell's original idea, one may for instance take into account absorption or diffusion phenomena at the boundary, thus obtaining
\begin{equation}\label{GWBC}
\left(\Delta u\right)_{|_{\partial\Omega}}  = -\frac{\partial u}{\partial \nu}+u_{|_{\partial\Omega}}+\Delta_{\partial \Omega} u_{|_{\partial\Omega}}\qquad \hbox{on }\partial \Omega
\end{equation}
(yet further nonlocal terms may be added, as we will see). Here we are using the notational convention that $\left(\Delta u\right)_{|_{\partial\Omega}}$ represents the trace of $\Delta u:\Omega\to \mathbb C$ at the boundary, whereas $\Delta_{\partial \Omega} u_{|_{\partial\Omega}}$ is obtained applying the Laplace--Beltrami operator (defined on the manifold $\partial \Omega$) to the trace of $u$. More recently, the Laplacian equipped with this more general class of boundary conditions has been studied in the $L^p$-framework, too, and a generation result appears in~\cite{CocFavGol08}.

But what to do if a variational structure is not apparent? One of the main ideas in~\cite{FavGolGol02,CocFavGol08} is that even when~\eqref{GWBC} is replaced by
\begin{equation}\label{eq:gwbcgen}
\Delta u  = -\beta \frac{\partial u}{\partial \nu}+\gamma u_{|_{\partial\Omega}}+\delta \Delta_{\partial \Omega} u_{|_{\partial\Omega}}\qquad \hbox{on }\partial \Omega
\end{equation}
for some functions $\beta,\gamma,\delta:\partial \Omega\to \mathbb C$, a smart renorming is sufficient to define a Lyapunov functional of the system and hence to show that the parabolic equation is governed by a dissipative semigroup:  This is done by endowing the naive state space $L^2(\Omega)\times L^2(\partial\Omega)$ by the inner product
\begin{equation}\label{eq:renorming}
(f|g):=\int_\Omega f_1 \overline{g_1} \ dx + \int_{\partial \Omega} f_2 \overline{g_2}\frac{d\sigma}{\beta}.
\end{equation} 
However, the introduction of the inner product in~\eqref{eq:renorming} is only feasible whenever $\frac{1}{\beta}$ is a well-behaved density, and in particular only if $\beta$ is positive. But even if we give up our wish for energy dissipation in the system, is the heat equation with boundary conditions~\eqref{eq:gwbcgen} still well-posed if no assumption is made on the sign of $\beta$? This issue has been studied in detail in~\cite{VazVit07}, and the answer is: \emph{Generally speaking, no}.

This result seems to be difficult to explain by variational methods. Instead, when it comes to discuss operators that are not self-adjoint -- or whose self-adjoint nature is not clear -- the techniques of Engel prove rather effective. Indeed, with Engel's methods checking the generator property for an operator matrix is \emph{morally} reduced to showing that two operators on $L^2(\Omega)$ and $L^2(\partial \Omega)$, respectively, are generators: This is usually easily done, since typically these operators are lower order perturbations of the Laplacian with Dirichlet boundary condition, and of the Dirichlet-to-Neumann operators, respectively. By ``morally'' we mean that in the previous literature one had to either restrict to the case $d=1$ (\cite{Mug06}) or to work in different functional spaces to avoid certain technical problems (\cite{EngFra05}). Our aim is to present how Engel's theory can be applied to show generation in an $L^p$-setting also for $d>1$. 

We begin by reviewing some known results in the theory of operator matrices with diagonal and non-diagonal domain in Sections~\ref{sec:diago} and~\ref{sec:nondiago}, respectively: We refer the interested reader to the monographs~\cite{Eng98,Tre08} for more thorough treatments of different aspects of this theory. We are going to study the long-time behavior of semigroups by applying suitable estimates on convolutions of vector-valued functions. As an application we present in Section~\ref{sec:3apps} some generalizations of known well- and ill-posedness results obtained in~\cite{BanBelRei06,VazVit07,CocFavGol08,VazVit11} and~\cite{Lum80a,Bel91} for scalar-valued functions on general domains and for diffusion equations on metric graphs, respectively.

\medskip
{\bf Notation.} Throughout the paper we denote by $[D(T)]$ the Banach space obtained by endowing the domain of a closed operator $T$ on a Banach space $W$ by its graph norm
\[
\|x\|_{[D(T)]}:=\|Tx\|_W + \|x\|_W.
\]
For two Banach spaces $W,Z$ and $\alpha \in [0,1]$ we denote by $[W,Z]_\alpha$ the complex interpolation space of order $\alpha$ between $W,Z$, cf.~\cite{LioMag72}. We denote by ${\mathcal L}_s(W,Z)$ the vector space of all bounded linear operators from $W$ to $Z$ endowed with the strong topology, whereas ${\mathcal L}(W,Z)$ denotes as usual the Banach space obtained endowing the same vector space with the norm topology.

\section{Operator matrices with diagonal domain}\label{sec:diago}

Whenever an abstract problem that is related to coupled systems of linear PDEs is studied, one commonly has to check whether an operator matrix \begin{equation}\label{opma}
\mathbf{A}:=\begin{pmatrix}
A & B\\
C & D
\end{pmatrix}
\end{equation}
generates a $C_0$-semigroup on a suitable product Banach space $E\times F$. The compact notation in~\eqref{opma} is meant to signify that $\mathbf A$ acts on elements of $E\times F$ by 
\[
{\mathbf A}:\begin{pmatrix}f_1\\ f_2 \end{pmatrix}\mapsto \begin{pmatrix}Af_1 + Bf_2\\ Cf_1 + Df_2\end{pmatrix}.
\]
Hence, $\mathbf A$ is itself an operator on $E\times F$ with domain
\[
D({\mathbf A})=\big(D(A)\cap D(C)\big)\times \big(D(D)\cap D(B)\big).
\]
The standing assumptions ob the linear operators $A,B,C,D$ are summarized in the following.

\begin{assums}\label{elem} Let us impose the following throughout this section.
\begin{enumerate}
\item $E$ and $F$ are Banach spaces.
\item ${A}:D({A})\subset E\to E$ is closed.
\item $D:D(D)\subset F\to F$ is closed.
\item $B:D(B)\subset F\to E$ and $C:D(C)\subset E\to F$.
\end{enumerate}
\end{assums}

In particular, by the above assumptions the domains of $A,D$ become Banach spaces with respect to the respective graph norm.

If $B,C$ are ``more unbounded'' than $A,D$, then the abstract Cauchy problem
\[
\frac{d {\mathbf u}}{dt}(t)={\mathbf Au}(t),\quad t\ge 0,\qquad {\mathbf u}(0)={\mathbf u}_0,
\]
has generally hyperbolic character. We will however rather focus on the opposite case and thus impose the following.

\begin{assum}\label{assum:lastmin}
$C$ is bounded from $[D(A)]$ to $\partial X$ and $B$ is bounded from $[D(D)]$ to $X$.
\end{assum}

If $B=0$ and $C=0$ it is an elementary exercise to check that $A$ and $D$ generate a semigroup on $E$ and $F$, respectively, if and only if $\mathbf A$ generates a semigroup on  ${E}\times {F}$. The following summarizes and slightly extends some results obtained in~\cite[\S~3]{Nag89}.

\begin{theo}\label{Nag89}
Under the Assumptions~\ref{elem} and~\ref{assum:lastmin} the following assertions hold for the operator matrix $\mathbf A$ defined in~\eqref{opma} with diagonal domain $D(\mathbf{A}):=D({A})\times D(D)$ on the product space ${\mathbf E}:=E\times F$.
\begin{enumerate}
\item Let 
\begin{itemize}
\item $B\in{\mathcal L}([D(D)],[D(A)])$, or else $B\in{\mathcal L}(F,E)$, and moreover
\item $C\in{\mathcal L}([D(A)],[D(D)])$, or else $C\in{\mathcal L}(E,F)$. 
\end{itemize}
If $A$ and $D$ both generate $C_0$-semigroups $(e^{tA})_{t\geq 0}$ on $E$ and $(e^{tD})_{t\geq 0}$ on $F$, respectively, then $\mathbf A$ generates a $C_0$-semigroup $(e^{t\mathbf A})_{t\geq 0}$ on ${\mathbf E}$.
\item Let 
\begin{itemize}
\item $B\in{\mathcal L}(F,E)$ or $B\in{\mathcal L}([D(D)],[D(A)])$, and additionally $C\in{\mathcal L}([D(A)],F)$; or else 
\item $C\in{\mathcal L}(E,F)$ or $C\in{\mathcal L}([D(A)],[D(D)])$, and additionally $B\in{\mathcal L}([D(D)],E)$.
\end{itemize}
Then $A$ and $D$ both generate analytic semigroups $(e^{tA})_{t\geq 0}$ on $E$ and $(e^{tD})_{t\geq 0}$ on $F$, respectively, if and only if $\mathbf A$ generates an analytic semigroup $(e^{t\mathbf A})_{t\geq 0}$ on ${\mathbf E}$. 
\item Let $A$ and $D$ both generate analytic semigroups $(e^{tA})_{t\geq 0}$ on $E$ and $(e^{tD})_{t\geq 0}$ on $F$, respectively. If there exists $\alpha\in (0,1)$ such that
\begin{itemize}
\item $B\in{\mathcal L}([D(D)],\big[[D(A)],E\big]_\alpha)$  and
\item $C\in{\mathcal L}([D(A)],\big[[D(D)],F\big]_\alpha)$,
\end{itemize}
then $\mathbf A$ generates an analytic semigroup $(e^{t\mathbf A})_{t\geq 0}$. Conversely, if $\mathbf A$ generates an analytic semigroup on $\mathbf E$ and 
$$\begin{pmatrix}0 & B\\ C & 0\end{pmatrix}\in{\mathcal L}\left([D({\mathbf A})],\big[[D({\mathbf A})],{\mathbf E}\big]_\alpha \right)$$
for some $\alpha\in(0,1)$, then also $A$ and $D$ generate semigroups on $E$ and $F$, respectively.

\item Assume 
\begin{itemize}
\item $B$ to be compact as an operator from $[D(A)]$ to $F$ and
\item $C$ to be compact as an operator from $[D(D)]$ to $E$.
\end{itemize}
Then $A$ and $D$ both generate analytic semigroups $(e^{tA})_{t\geq 0}$ on $E$ and $(e^{tD})_{t\geq 0}$ on $F$, respectively, if and only if $\mathbf A$ generates an analytic semigroup $(e^{t\mathbf A})_{t\geq 0}$ on ${\mathbf E}$.
\end{enumerate}
If any of the above assertions hold with $B=0$ or $C=0$, then 
\begin{equation}\label{formsem}
e^{t\mathbf A}=
\begin{pmatrix}
e^{t{A}} & 0\\
{R}(t) & e^{tD}
\end{pmatrix}\quad \hbox{or}\quad
e^{t\mathbf A}=
\begin{pmatrix}
e^{t{A}} & S(t)\\
0 & e^{tD}
\end{pmatrix},\qquad t\geq 0,
\end{equation}
respectively, where the operator families $(R(t))_{t\ge 0},(S(t))_{t\ge 0}$ (necessarily consisting of bounded operators and strongly continuous) are defined by
\begin{equation}\label{rst}
{R}(t):=\int_0^t e^{(t-s)D}{C}e^{sA} ds,\qquad S(t):=\int_0^t e^{(t-s)A}B  e^{sD}ds,\qquad t\geq 0.
\end{equation}
\end{theo}

Most of these assertions are proved applying standard perturbation theorems to the case of operator matrices. The above statements do not exhaust at all the possible applications of classical perturbation theory: Further results may be applied, like those discussed in~\cite[Chapter~3]{EngNag00} or~\cite[Chapters 4--5]{BanArl06}, but their applications is usually more delicate and often requires detailed knowledge of the semigroups $(e^{tA})_{t\ge 0}$ and $(e^{tD})_{t\ge 0}$ -- detailed to an extent that is seldom given in concrete applications.

\begin{proof}
The assertion (2) and the identity~\eqref{formsem} are~\cite[Prop.~3.1 and Cor.~3.2]{Nag89}. In order to prove (1) and (3), write
\begin{equation}\label{first}
{\mathbf A}={\mathbf A}_0+{\mathbf A}_1:=\begin{pmatrix} A & 0\\ 0 & D\end{pmatrix} + \begin{pmatrix} 0 & B\\ C & 0\end{pmatrix}
\end{equation}
and observe that the first addend ${\mathbf A}_0$ on the right hand side has diagonal domain $D(A)\times D(C)$ and generates a semigroup. 

If this semigroup is merely strongly continuous, like in (1), then the generator property follows applying the Bounded Perturbation Theorem and/or perturbation result obtained by  Desch and Schappacher in~\cite[pag.~335]{DesSch84}, cf.\ also~\cite[Thm.~III.1.3 and Cor.~III.1.5]{EngNag00}.

If on the other hand the semigroup generated by ${\mathbf A}_0$ is analytic, one can easily describe the complex interpolation spaces between
$[D({\mathbf A}_0)]$ and $\mathbf E$ by
\[
\big[[D({\mathbf A}_0)],{\mathbf E}\big]_\alpha=\big[[D(A)],E\big]_\alpha \times \big[[D(D)],F\big]_\alpha,\qquad \alpha\in (0,1).
\]
Now, under the assumptions of (3)-(4) ${\mathbf A}_1$ is  bounded from $[D({\mathbf A})]$ to $\big[[D({\mathbf A})],{\mathbf E}\big]_\alpha$, or else compact  from $[D({\mathbf A})]$ to $\mathbf X$, and the claim follows by further perturbation results obtained by Desch and Schappacher in~\cite[pag.~338]{DesSch84} and~\cite{DesSch88}, cf.\ also~\cite[Thm.~3.7.25]{AreBatHie01}.
\end{proof}

\begin{rem}\label{bpt} Let the operators $B,C$ be bounded.

1) By the Bounded Perturbation Theorem one obtains that if $M_A,M_D,\epsilon_A,\epsilon_D$ are constants such that 
\[
\Vert e^{tA}\Vert_{\mathcal L(E)} \leq M_A e^{\epsilon_A t}\quad\hbox{ and }\Vert e^{tD}\Vert_{\mathcal L(F)} \leq M_D e^{\epsilon_D t},\qquad t\geq 0,
\]
then 
\begin{equation}\label{eq:altab}
\Vert e^{t\mathbf A}\Vert_{\mathcal L(\mathbf E)} \leq M e^{(\epsilon+M\max\{\Vert B\Vert,\Vert C\Vert\})t},\qquad t\geq 0,
\end{equation}
where
\[
M:=\max\{M_A,M_D\}\qquad \hbox{and}\qquad \epsilon:=\max\{\epsilon_A,\epsilon_D\}.
\]
In particular, $(e^{t{\mathbf A}})_{t\geq 0}$ is uniformly exponentially stable -- i.e., 
\[
\|e^{t{\mathbf A}}\|_{\mathcal L(\mathbf E)}\le Me^{\omega t},\qquad t\ge 0,
\] 
for some $\omega<0$ -- 
provided that
\[
M\max\{\Vert B\Vert,\Vert C\Vert\}<-\epsilon.
\]

2) By the Datko--Pazy Theorem (cf.~\cite[Thm.~V.1.8]{EngNag00}) a $C_0$-semigroup on a Banach space $W$ is uniformly exponentially stable if and only if it is of class $L^p({\mathbb R}_+,{\mathcal L}_s({W}))$ for some/all $p\in [1,\infty)$: In this case, its generator is invertible by~\cite[Prop.~IV.2.2]{EngNag00}. This has been used in~\cite{Mug11} to show that if both $\epsilon_A,\epsilon_B$ are strictly negative (i.e., if $(e^{tA})_{t\ge 0},(e^{tD})_{t\ge 0}$ are both uniformly exponentially stable) and
\begin{equation}\label{stab}
M_A M_D \Vert B\Vert \Vert C\Vert< {\epsilon_A\epsilon_B},
\end{equation}
then $( e^{t\mathbf A})_{t\ge 0}$ is uniformly exponentially stable, too. In comparison with the criterion in (1), this is particularly useful whenever $M_A\not=M_D$ or $\epsilon_A\neq \epsilon_D$. If in particular $A$ and $D$ are negative definite operators with compact resolvent whose largest (strictly negative) eigenvalues are $\lambda_A$ and $\lambda_D$, respectively, then the above criterion simply says that $( e^{t\mathbf A})_{t\ge 0}$ is uniformly exponentially stable provided that 
\begin{equation*}\label{stab-}
\Vert B\Vert \Vert C\Vert< {\lambda_A\lambda_D}.
\end{equation*}
\end{rem}

The formulae in~\eqref{formsem} are useful to deduce further information on $(e^{t\mathbf A})_{t\ge 0}$. We will also devote our attention to asymptotic almost periodicity, see e.g.~\cite[\S~4.7 and~\S~5.4]{AreBatHie01} for an explanation and discussion of this notion and its relevance for both bounded uniformly continuous functions and bounded operators semigroups. Roughly speaking, a function on $\mathbb R_+$ with values in a Banach space is asymptotically almost periodic (shortly: AAP) if it is the direct sum of a continuous function vanishing at $+\infty$ and a limit of linear combinations of rotations. In particular, a uniformly exponentially stable semigroup is trivially AAP, too. By~\cite[Thm.~4.7.4]{AreBatHie01}, each bounded compact semigroup, and in particular each analytic bounded semigroup whose generator has compact resolvent, is AAP.

If the operator matrix $\mathbf A$ is upper or lower triangular, the form of $(R(t))_{t\geq 0}$ and $(S(t))_{t\geq 0}$ allows us to apply known results on convolutions of operator valued mappings.

\begin{theo}\label{boundNag89}
Assume any of the cases in Theorem~\ref{Nag89} to hold with $B=0$, and in particular $\mathbf A$ and hence $A,D$ to be semigroup generators.

\begin{enumerate}[(1)]
\item If $(e^{tA})_{t\ge 0},(e^{tD})_{t\ge 0}$ are both uniformly exponentially stable, then so is $(e^{t\bf A})_{t\ge 0}$.
\item Let $(e^{tA})_{t\geq0}$ and $(e^{tD})_{t\geq 0}$ be bounded and uniformly exponentially stable, respectively.  Then the following hold.
\begin{enumerate}
\item The semigroup $(e^{t\mathbf A})_{t\geq0}$ is bounded.
\item If $(e^{t A})_{t\ge 0}$ is AAP, then so is $(e^{t\mathbf A})_{t\ge 0}$.
\item If $\lim\limits_{t\to\infty}e^{t{A}}$ exists (resp., exists and is equal 0) in the strong operator topology, then so does $\lim\limits_{t\to\infty}R(t)$ and
\[
\lim_{t\to\infty} R(t)x=D^{-1}C\lim_{t\to\infty}e^{t{A}}x,\qquad x\in E.
\]
\item If $(e^{t{A}})_{t\geq0}$ is uniformly exponentially stable, then so is $(e^{t\mathbf A})_{t\geq0}$.
\end{enumerate}
\item Let $(e^{t{A}})_{t\geq0}$ and $(e^{tD})_{t\geq 0}$ be uniformly exponentially stable and bounded, respectively.  Then the following hold.
\begin{enumerate}
\item The semigroup $(e^{t{\mathbf A}})_{t\geq0}$ is bounded.
\item If $(e^{tD})_{t\geq 0}$ is AAP, then $(R(t))_{t\geq0}$ is AAP.
\item If $\lim\limits_{t\to\infty}e^{tD}$ exists (resp., exists and is equal 0) in the strong operator topology, then so does $\lim\limits_{t\to\infty}e^{t\mathbf A}$.
\end{enumerate}
\item Assume that the spectra of $A,D$ satisfy
\begin{equation}\label{eq:nonres-0}
i{\mathbb R}\cap \sigma(A)\cap \sigma (D)=\emptyset .
\end{equation}
If moreover $(e^{tA})_{t\geq 0}$ and $(e^{tD})_{t\geq0}$ are bounded, then the following assertions hold.
\begin{enumerate}
\item If $(e^{tD})_{t\geq0}$ is analytic, then the orbit $(R(t))_{t\geq0}$ is bounded.
\item Let $(R(t))_{t\geq0}$ be bounded. If $(e^{tD})_{t\geq0}$ and $(e^{tA})_{t\geq 0}$ are AAP, then so is $(e^{t\mathbf A})_{t\geq0}$.
\item Let $(R(t))_{t\geq0}$ be bounded. If $\lim\limits_{t\to\infty}e^{tD}$ and $\lim\limits_{t\to\infty}e^{tA}$ exist (resp., exist and are equal 0) in the strong operator topology, then so does $\lim\limits_{t\to\infty}e^{t\mathbf A}$.
\end{enumerate}

\end{enumerate}
\end{theo}

Of course, analogous results hold whenever $C=0$ and $B\in{\mathcal L}(F,E)$. 

\begin{proof}
(1) Exponential stability of $(e^{t\bf A})_{t\ge 0}$ follows directly from  Remark~\ref{bpt}.(2).

Let now $x\in E$. By~\eqref{rst},
\[
R(t)x=(e^{\cdot D}*f_x)(t),\qquad t>0,
\]
where 
\begin{equation}\label{def:f}
f_x:=Ce^{\cdot {A}}x,\qquad x\in E.
\end{equation}
Observe that in fact $f_x$ is well-defined by Theorem~\ref{Nag89} and indeed $f_x\in L^1_{loc}({\mathbb R}_+,F)$. Now, the assertions in (2) and (3) follow from~\cite[Prop.~5.6.1]{AreBatHie01} and~\cite[Prop.~5.6.4]{AreBatHie01}, respectively.

(4) For any $x\in E$ the Laplace transform  $\hat{f}_x$ of $(f_x(t))_{t\ge 0}$ in~\eqref{def:f} is  $CR(\lambda, A)x$, ${\rm Re}\lambda>0$: Thus the \emph{half-line spectrum} ${\rm sp}(f_x)$ of $f_x$, defined as in~\cite[\S~4.4]{AreBatHie01}, is $\{\eta\in {\mathbb R}: i\eta\in\sigma(A)\}$. Then the claims follow from~\cite[Thm.~5.6.5 and Thm.~5.6.6]{AreBatHie01}.
\end{proof}

\begin{rem}\label{findings}
We can summarize our findings about the semigroup $(e^{t\mathbf A})_{t\ge 0}$ (for $B=0$ and $C\in \mathcal L([D(A)],F)$ in tabular form as follows:\\

\begin{center}
	\begin{tabular}{r|l|l|l}
		\backslashbox{$e^{tA}$}{$e^{tD}$} & unif.\ exp.\ stable & bounded & AAP \\ \hline
unif.\ exp.\ stable & unif.\ exp.\ stable & bounded & AAP  \\ \hline
		bounded & bounded &  bounded* &  ?\\ \hline
		AAP & AAP  & ?  &  AAP*
			\end{tabular}
	\label{tab:resume}
\end{center}

\smallskip
(* holds if additionally $(e^{tD})_{t\ge 0}$ is analytic and~\eqref{eq:nonres-0} is satisfied.)

If instead $C=0$ and $B\in \mathcal L([D(D)],E)$, then the same table prevails, but * holds instead provided $(e^{tA})_{t\ge 0}$ is analytic and~\eqref{eq:nonres-0} is satisfied.
\end{rem}

\section{Operator matrices with non-diagonal domain}\label{sec:nondiago}

In this section we are going to deduce results similar to those of Section~\ref{sec:diago} for the same operator matrix 
$${\mathbf A}:=\begin{pmatrix}
A & 0\\
C & D	
\end{pmatrix} ,$$
defined however on a smaller, \emph{coupled} domain 
\begin{equation}\label{opmadom}
D({\mathbf A}):=\left\{
\begin{pmatrix}
{u}\\{x}
\end{pmatrix}
\in \big(D({A})\cap D(C)\big)\times \big(D(D)\cap D(B)\big):\ Lu=x\right\},
\end{equation}
where $L$ is an operator that, in view of applications, we think of as a \emph{boundary operator} (e.g., a normal derivative operator or a point evaluation functional). Throughout this section we consistently replace the assumptions imposed in Section~\ref{sec:diago} by the following ones.

\begin{assums}\label{basic} $ $
\begin{enumerate}[(1)]
\item $X,\partial X$ are Banach spaces.
\item ${A}:D({A})\subset {X}\to {X}$.
\item $L:D({A})\subset X\to \partial X$ is surjective.
\item $A_0:={A}_{|\ker \; L}$ has nonempty resolvent set $\rho(A_0)$.
\item $\begin{pmatrix}{A}\\{L}\end{pmatrix}:D(A)\to {X}\times {\partial X}$ is closed.
\item $D:D(D)\subset \partial X\to \partial X$ is closed.
\item $B:D(B)\subset \partial X\to X$ and $C:D(C)\subset X\to \partial X$.
\item $C$ is bounded from $[D(A_0)]$ to $\partial X$.
\end{enumerate}
\end{assums}

We denote by $[D(A)_L]$ the Banach space obtained by endowing $D(A)$ with the graph norm of $A\choose L$. Clearly, $[D(A_0)]\hookrightarrow [D(A)_L]$. Under the Assumptions~\ref{basic} $A_0$ is closed, and furthermore $[D(A_0)]\hookrightarrow [D(A)_L ]$. Accordingly, $\mathbf A$ is closed and $[D({\mathbf A})]\hookrightarrow [D(A)_L ]\times \left([D(D)]\cap \partial Y\right)$.

\begin{exa}\label{exa-00}
The most relevant application of our results will be to the Laplacian with Wentzell-type dynamic boundary conditions. We are going to discussed this in detail in Section~\ref{wbc}. However, let us already point out that the above assumptions are satisfied for the following choice of spaces and operators, whenever $\Omega$ is an open domain of $\mathbb R^d$ whose boundary we assume to be nonempty (to avoid trivialities) and smooth:
\[
X:=L^2(\Omega),\qquad \partial X:=L^2(\partial \Omega),
\]
as well as
\[
A:=\Delta,\quad D(A):= \left\{ u\in H^\frac12(\Omega):\Delta u\in L^2(\Omega)\right\},\qquad B=0.
\]
It is known that the trace operator
\[
L:u\mapsto u_{|_{\partial\Omega}}
\]
is bounded and surjective from $D(A)$ to $\partial X$ whenever $\partial\Omega$ is smooth enough, cf.~\cite[Vol.~I, Thm.~2.7.4]{LioMag72}. Moreover, the closedness of ${A\choose L}$ holds by interior estimates for general elliptic operators (a short proof of this can be found in \cite[\S~3]{CasEngNag03}). We see that $A_0$ is a (weakly defined) realization of the Laplacian with Dirichlet boundary conditions, the generator of an analytic semigroup on $L^2(\Omega)$: It is invertible (i.e., $0\in \rho(A_0)$) and its domain is
\[
D(A_0):=H^2(\Omega)\cap H^1_0(\Omega)
\]
We finally introduce the normal derivative operator
\[
C:u\mapsto \frac{\partial u}{\partial \nu},\qquad D(C):= \left\{ u\in D(A):\frac{\partial u}{\partial \nu}\in L^2(\partial \Omega)\right\}.
\]
Again by the results of~\cite{LioMag72}, $D(C)\subset H^{\frac{3}{2}}(\Omega)$.
\end{exa}

Let us consider the abstract elliptic boundary value problem
\begin{equation}\tag{AEP}
\left\{
\begin{array}{rll}
{Au}&=&\lambda {u},\\
{Lu}&=&{x}.
\end{array}
\right.
\end{equation}
The following is a slight modification of a result due to Greiner, cf.~\cite[Lemma~2.3]{CasEngNag03}.

\begin{lemma}\label{dirichdef}
Under the Assumptions~\ref{basic}, let $\lambda\in\rho(A_0)$. Then the problem $({\rm AEP})$ admits a unique solution ${u}:=D^{{A},{L}}_\lambda {x}$ for all ${x}\in {\partial X}$. Moreover, the operator $D^{A,L}_\lambda$ is bounded from $\partial X$ to $W$ for all Banach spaces $W$ that satisfy $D(A^k)\subset W\hookrightarrow {X}$ for some $k\in \mathbb N$. 
\end{lemma}

\begin{assums}\label{basic-2} $ $
\begin{enumerate}[(1)]
\item $\partial Y$ is a Banach space that is continuously and densely embedded in $\partial X$.
\item $L_{|D(C)}:D(C)\to \partial Y$ is surjective.
\item $\begin{pmatrix}{A}\\{L}\end{pmatrix}_{|D(C)}:D(C)\to {X}\times {\partial Y}$ is closed.
\end{enumerate}
\end{assums}

Applying Lemma~\ref{dirichdef} to $\partial Y$ instead of $\partial X$, Assumptions~\ref{basic-2} ensure that $D^{A,L}_\lambda$ has a well-behaved restriction from $\partial Y$ to $Z$ for all Banach spaces $Z$ that satisfy $D(C)\subset {Z}\hookrightarrow {X}$. Furthermore, $D(D)\cap \partial Y$ becomes a Banach space for the natural sum norm, and so does $D(C)$ with respect to
\[
\|u\|_{[D(C)_{A,L}]}:=\|u\|_X + \|Au\|_X+\|Cu\|_{\partial X}+\|Lu\|_{\partial Y}.
\]
Hence, for all $\lambda \in \rho(A_0)$ $D^{A,L}_\lambda$ is well-defined and bounded also from ${\partial Y}$ to $[D(C)_{A,L}]$ and hence $CD^{A,L}_\lambda$ is bounded from ${\partial Y}$ to $\partial X$.

If now $\lambda\in\rho(A_0)$, then $D+{C}D^{{A},{L}}_\lambda$ with domain $D(D)\cap \partial Y$ is well-defined and the factorization
\begin{equation}\label{klen1}
{\mathbf{A}}-\lambda
= {\mathbf{A}}_\lambda {\mathbf M}_\lambda\\
:=\begin{pmatrix}
A_0-\lambda & B\\
C & D+{C}D^{{A},{L}}_\lambda-\lambda
\end{pmatrix}
\begin{pmatrix}
\Id & -D^{A,L}_\lambda\\
0 & \Id
\end{pmatrix}
\end{equation}
is easily seen to hold -- this observation goes back to~\cite{Eng97,Eng99,KraMugNag03b}. By Lemma~\ref{dirichdef} the operator ${\mathbf M}_\lambda$ is is an isomorphism on ${\mathbf X}=X\times \partial X$ for all $\lambda\in \rho(A_0)$ with 
\[
{\mathbf M}_\lambda^{-1}=\begin{pmatrix}
\Id & D^{A,L}_\lambda\\
0 & \Id
\end{pmatrix}.
\]

\begin{exa}\label{exa-01}
Let $\lambda\in \rho(A_0)$. In the setting of Example~\ref{exa-00} $CD^{A,L}_\lambda$ agrees by definition with the so-called \emph{Dirichlet-to-Neumann operator} $\DN_\lambda$ -- i.e, the pseudo-differential operator defined by 
\[
\DN_\lambda f:=-\frac{\partial u}{\partial \nu},\qquad \hbox{where $u$ is the weak solution of }
\left\{
\begin{array}{rcll}
\Delta u&=&\lambda u \qquad &\hbox{in }\Omega,\\
u&=&f\qquad &\hbox{on }\partial \Omega.
\end{array}
\right.
\]
This operator occurs often in the contexts of PDEs and control theory: We are going to discuss some of its properties in more detail in Remark~\ref{rem:spectr} below. 

Remarkably, also the Dirichlet-to-Neumann operator $\DN_0$ is known to be associated with a Markov process, cf.~\cite[Thm.~9.1]{SatUen65}. This motivates to study the very general Wentzell-type boundary condition
\begin{equation}\tag{GWBC}
 \frac{\partial u}{\partial t}=-\beta\frac{\partial u}{\partial \nu}+\eta_1 \Delta_{\partial \Omega}u_{|_{\partial\Omega}}+\eta_2 \DN_0u_{|_{\partial\Omega}}+\eta_3  u_{|_{\partial\Omega}}\qquad \hbox{on }\partial \Omega,
\end{equation}
for $\beta,\eta_1,\eta_2,\eta_3\in \mathbb R$: This is e.g.\ in analogy with the setting of~\cite[Exa.~5.9]{Pos13}, once one remembers that the Dirichlet-to-Neumann operator agrees with $-(-\Delta_{\partial \Omega})^\frac12$ up to a lower order perturbation, cf.~\cite[Prop.~C.1, pag.~453]{Tay96}). Let us also emphasize that the seemingly related boundary condition
\[
\frac{\partial u}{\partial \nu}+\DN_0 u_{|_{\partial\Omega}}=0
\]
is well-known in the theory of self-adjoint extensions: Indeed, it defines the so-called \emph{Krein--von Neumann extension} of the Laplacian, cf.~\cite{AshGesMit13}. We are therefore led to consider the operator
\[
D:=\eta_1 \Delta_{\partial \Omega}+\eta_2 \DN_0+\eta_3 \Id\ ,
\]
which is closed on $L^2(\partial \Omega)$
\begin{itemize}
\item with domain $H^2(\partial \Omega)$ if $\eta_1\neq 0$;
\item with domain $H^1(\partial \Omega)$ if $\eta_1= 0$ but $\eta_2\neq 0$;
\item with domain $L^2(\partial \Omega)$ if $\eta_1= \eta_2=0$.
\end{itemize}
\end{exa}

In view of~\eqref{klen1} we promptly deduce the following.

\begin{cor}\label{factormain2}
Under the Assumptions~\ref{basic} and~\ref{basic-2}, for all $\lambda\in\rho(A_0)$ the operator matrix ${\mathbf A}-\lambda$ is similar to
\begin{equation}\label{simileq}
{\mathbf M}_\lambda {\mathbf A}_\lambda =\tilde{\mathbf A}_\lambda:=
\begin{pmatrix}
A_0-{D^{A,L}_\lambda C}-\lambda & B-D^{A,L}_\lambda(D+ C D^{A,L}_\lambda -\lambda)\\
{C} & D+ C D^{A,L}_\lambda -\lambda
\end{pmatrix},
\end{equation}
with domain
$$D(\tilde{\mathbf A}_\lambda):=D(A_0)\times (D(D)\cap \partial Y).$$
In particular, if $\mathbf A$ has nonempty resolvent set, then $\mathbf A$ has compact resolvent if and only if the embeddings $[D(A_0)]\hookrightarrow {X}$ and $[D(D)]\cap \partial Y\hookrightarrow {\partial X}$ are both compact.
\end{cor}

By Corollary~\ref{factormain2} the operator matrix with coupled domain $\bf A$ is a generator on $\bf X$ if and only if the similar operator matrix $\tilde{\bf A}_\lambda$ -- which has diagonal domain! -- is a generator on the same space. The main idea is now to regard $\tilde{\bf A}_\lambda$ as the sum
\begin{equation}\label{eq:maindec}
\begin{array}{l}
\tilde{\mathbf A}_\lambda=
\begin{pmatrix}
A_0 & 0\\
C & D+{C}D^{{A},{L}}_\lambda\\
\end{pmatrix} +\begin{pmatrix}
-D^{A,L}_\lambda C-\lambda & B-{D^{{A},{L}}_\lambda}(D+{C}D^{{A},{L}}_\lambda+\lambda D^{{A},{L}}_\lambda)\\
0 & -\lambda
\end{pmatrix} :
\end{array}
\end{equation}
Now, several sets of conditions implying that $\tilde{\bf A}_\lambda$ is a generator become available as a consequence of Theorem~\ref{Nag89}. We only explicitly mention those two that prove most relevant in view of our applications.

\begin{theo}\label{gener} 
Under the Assumptions~\ref{basic}-\ref{basic-2} the following assertions hold.
\begin{enumerate}
\item Let $A_0$ and $D+{C}D^{{A},{L}}_\lambda$ generate analytic semigroups on $X$ and $\partial X$, respectively, for some $\lambda\in \rho(A_0)$. Let for some $0<\alpha<1$ the complex interpolation spaces associated with $A_0$ and $D+{C}D^{{A},{L}}_\lambda$ satisfy
\begin{itemize}
\item $D^{A,L}_\lambda(\partial X)\hookrightarrow \big[[D(A_0)],{X}\big]_\alpha$,
\item $B\in{\mathcal L}\left([D(D)\cap\partial Y],\big[[D(A_0)],X\big]_\alpha\right)$.
\end{itemize}
Then $\mathbf A$ generates an analytic semigroup on $\mathbf X$. 

Conversely, let $\tilde{\mathbf A}_\lambda$ generates an analytic semigroup -- say, with diagonal domain $W\times Z$, for two subsets $W$ of $X$ and $Z$ of $\partial X$ -- and assume that
\begin{itemize}
\item $D^{A,L}_\lambda(\partial X)\hookrightarrow \big[W,X\big]_\alpha$,
\item $B\in{\mathcal L}\left(Z,\big[[D(A_0)],X\big]_\alpha\right)$.
\end{itemize}
Then also $A_0$ and $D+{C}D^{{A},{L}}_\lambda$ generate an analytic semigroup on $X,\partial X$, respectively. 
\item Assume for some $\lambda\in \rho(A_0)$
\begin{itemize}
\item $D^{A,L}_\lambda$ to be compact as an operator from $\partial X$ to $X$,
\item $B$ to be compact as an operator from $[D(D)]\cap\partial Y$ to $X$.
\end{itemize}
Then $A_0$ and $D+{C}D^{{A},{L}}_\lambda$ generate analytic semigroups on $X$ and $\partial X$, respectively, if and only if $\mathbf A$ generates an analytic semigroup on $\mathbf X$.
\end{enumerate}
\end{theo} 

\begin{proof} 
Take $\lambda\in\rho(A_0)$. We decompose $\tilde{\mathbf A}_\lambda$ as in~\eqref{eq:maindec}. By assumption $C$ is bounded from $[D(A_0)]$ to $\partial X$, hence by Theorem~\ref{Nag89}.(3)-(4) the first addend on the right-hand side of~\eqref{eq:maindec} generates an analytic semigroup on $\mathbf X$ if and only if its diagonal entries generate an analytic semigroup on $X,\partial X$, respectively. 

It remains to observe that the second addend on the right hand side is a well-behaved perturbation under either set of assumptions. Indeed, if (1) holds, then for $\alpha\in(0,1)$ the complex interpolation space corresponding to the leading term is $\big[[D(\tilde{\mathbf A}_\lambda)], {\mathbf X}\big]_\alpha=[D(A_0),{X}]_\alpha\times \big[[D(D)]\cap \partial Y,\partial X\big]_\alpha$. Thus, by assumption the second addend on the right-hand side is bounded from $[D(\tilde{\mathbf A}_\lambda)]$ to $\big[[D(\tilde{\mathbf A}_\lambda)],{\mathbf X}\big]_\alpha$. If instead (2) holds, then the second addend on the right-hand side is compact from $[D(\tilde{\mathbf A}_\lambda)]$ to ${\mathbf X}$. In either case, the claim follows by the perturbation theorems due to Desch and Schappacher already exploited in the proof of Theorem~\ref{Nag89}.(3)-(4).
\end{proof}


\begin{rem}\label{exa:misc}
Let us assume that 
\[
A_0 \hbox{ is invertible}.
\]
If $C=0$, then we can allow for slightly more general operators $B$ than in Theorem~\ref{gener}. More precisely, let $A_0,D$ generate analytic semigroups and let
\[
B\hbox{ be bounded from }[D(D)]\cap \partial Y\hbox{ to }X\qquad \hbox{and}\qquad C=0.
\]
Then by Lemma~\ref{factormain2} the operator ${\mathbf A}$ is similar to
\begin{equation}\label{formgen}
\tilde{\mathbf A}_0:=
\begin{pmatrix}
A_0 & B-D_0^{A,L} D\\
0 & D
\end{pmatrix},
\end{equation}
with diagonal domain $D(\tilde{\mathbf A}_0):=D(A_0)\times D(D)$. By Theorem~\ref{Nag89} we conclude that $\mathbf A$ is a generator and in fact $e^{t\mathbf A}$ is for all $t\ge 0$ similar to
\begin{equation*}
\begin{pmatrix}
e^{t{A_0}} & S(t)\\
0 & e^{tD}
\end{pmatrix},\qquad t\geq 0,
\end{equation*}
where $(S(t))_{t\geq 0}$ is the strongly continuous family of convolution operators in~\eqref{rst}. Let us collect some direct consequences of the stability results obtained in Section~\ref{sec:diago}.

\begin{enumerate}
\item Uniform exponential stability of $(e^{t{A_0}})_{t\ge 0}$ and $(e^{t{D}})_{t\ge 0}$ is a necessary and sufficient condition for $(e^{t\mathbf A})_{t\ge 0}$ to be uniformly exponentially stable.
\item Let $(e^{t{A_0}})_{t\ge 0}$ be uniformly exponentially stable. If $(e^{tD})_{t\geq0}$ is bounded/AAP, then so is $(e^{t\mathbf A})_{t\geq0}$.
\item Let $(e^{t{A_0}})_{t\ge 0}$ be analytic and let
\begin{equation}\label{eq:nonres-bis}
i{\mathbb R}\cap \sigma(A_0)\cap \sigma (D)=\emptyset .
\end{equation}
If both $(e^{tA_0})_{t\geq0},(e^{tD})_{t\geq0}$ are bounded, then so is $(e^{t\mathbf A})_{t\geq0}$.
\end{enumerate}
The above assertions complement some results obtained in~\cite[\S~5]{CasEngNag03} by means of Perron--Frobenius-type results, and in particular under positivity assumptions: Observe that~\eqref{eq:nonres-bis} is satisfied in particular whenever $A_0$ or $D$ is self-adjoint and invertible.

We can also consider the case of a complete matrix $\mathbf A$ provided that 
\[
B \hbox{ is bounded from }\partial X\hbox{ to }X\qquad \hbox{and}\qquad C \hbox{ is bounded from }X\hbox{ to }\partial X.
\]
Then Theorem~\ref{gener} applies again and $\mathbf A$ generates an analytic semigroup if so do both $A_0$ (and hence $A_0+D_0^{A,L}C$) and $D+CD_0^{A,L}$.
Let  $M_A,M_B\geq 1$ and $\epsilon_A,\epsilon_B<0$ be constants such that
\begin{equation*}
\Vert e^{t(A_0-D^{A,L}_0C)}\Vert \leq M_Ae^{\epsilon_A t}\qquad \hbox{and}\qquad \Vert e^{t(D+CD^{A,L}_0)}\Vert \leq M_Be^{\epsilon_B t},\qquad t\geq 0.
\end{equation*}
Then $(e^{t{\mathbf A}})_{t\ge 0}$ is uniformly exponentially stable, too, provided that the estimate
\begin{equation}\label{stabbar}
\Vert C\Vert_{\mathcal L(X,\partial X)} \Vert B-D^{A,L}_0(D+CD_0^{A,L})\Vert_{\mathcal L(\partial X,X)}< \frac{\epsilon_A\epsilon_B}{M_A M_B}
\end{equation}
holds.

\end{rem}

\section{Applications}\label{sec:3apps}

\subsection{The Laplacian on a domain}\label{wbc}
The generation result of Section~\ref{sec:nondiago} can be applied in order to discuss the Laplacian on $\Omega\subset \mathbb R^n$, $n\ge 2$, equipped with the general Wentzell boundary conditions $\rm(GWBC)$ introduced in Example~\ref{exa-01}. The natural $L^2$-realization of the Laplacian equipped with (GWBC) is the operator matrix
$${\mathbf A}:=\begin{pmatrix}
\Delta & 0\\
-\beta\frac{\partial}{\partial\nu} & \eta_1 \Delta_{\partial \Omega}+\eta_2 \DN_0+\eta_3 \Id
\end{pmatrix}$$
with domain
\begin{eqnarray*}
&&D({\mathbf A}):=\left\{\begin{pmatrix}
u\\ x	
\end{pmatrix}
\in H^\frac12(\Omega)\times L^2(\partial\Omega): \Delta u\in L^2(\Omega),\right.\\
&&\qquad\qquad\qquad\left.
 \eta_1\Delta_{\partial \Omega} u_{|_{\partial\Omega}}\in L^2(\partial \Omega), \beta\frac{\partial u}{\partial\nu}\in L^2(\partial\Omega),\ \eta_2 \DN u_{|_{\partial\Omega}}\in L^2(\partial\Omega)\;\hbox{and } u_{|_{\partial\Omega}}=x\right\},
\end{eqnarray*}
for $\beta,\eta_1,\eta_2,\eta_3\in \mathbb R$. In order to apply the abstract results of the previous section, consider $\mathbf A$ as an operator matrix $\mathbf A$ with domain $D({\mathbf A})$ defined as in~\eqref{opma}--\eqref{opmadom}. Here we adopt the setting of Examples~\ref{exa-00}-\ref{exa-01}, and assume in particular $\partial Omega$ to be smooth.

Since $0$ lies in the resolvent set of $A_0$, by Lemma~\ref{dirichdef} the operator $D_0^{A,L}$ is a well-defined, bounded operator from $\partial X$ to $X$. 
By \cite[Vol.~II, (4.14.32)]{LioMag72}, we obtain that 
$$D^{A,L}_0(\partial X)\hookrightarrow H^\frac{1}{2}(\Omega)=\big[[D(A_0)],X\big]_{\frac{3}{4}}.$$

We still need to take a closer look at $CD^{A,L}_0=\beta\DN_0$: It was proved in~\cite[Lemma~6.2]{SatUen65} that 
$$\DN_0 :f\mapsto-\frac{\partial}{\partial \nu} D^{A,L}_0 f,$$
which has already been discussed in detail in Example~\ref{exa-01}, generates a semigroup on $C(\partial \Omega)$; while the \emph{analyticity} of this semigroup is the main result in~\cite{Esc94}. The $L^2$-realization of $\DN_0$ has been shown in~\cite[Thm.~4]{Kun70} (and independently in~\cite[\S~4]{EmaLaa06}) to generate an analytic semigroup on $\partial X=L^2(\Omega)$ if $\beta>0$; $\beta>0$ is also a necessary condition as long as $\Omega \subset \mathbb R^n$ with $n\ge 2$, since $\DN_0$ is then unbounded -- in fact, a pseudo-differential operator of order 1. (The case in which $\partial\Omega$ is Lipschitz is treated in~\cite{AreMaz12}.) We can thus complement some known ill-posedness results, cf.~\cite[Thm.~9]{BanBelRei06} and~\cite[Thm.~1]{VazVit08}, which deal with the case of $\eta_1=0$; while well-posedness in $H^1(\Omega)$ has been established in~\cite{VazVit11} in the case of $\eta_1>0$.

\begin{theo}\label{theo:main-ww}
Let $\Omega$ be an open, bounded domain of ${\mathbb R}^n$  ($n\ge 2$) with smooth boundary $\partial\Omega$. Let $\beta,\eta_1,\eta_2,\eta_3\in \mathbb R$. Then the operator matrix $\bf A$ generates an analytic (compact) semigroup on $L^2(\Omega)\times L^2(\partial \Omega)$ if and only if 
\[
CD^{A,L}_0+D:=\beta\DN_0+\eta_1 \Delta_{\partial \Omega}+\eta_2 \DN_0+\eta_3  \Id
\]
generates an analytic semigroup on $L^2(\partial \Omega)$: This is in turn the case if and only if either of the following conditions holds:
\begin{itemize}
\item $\eta_1> 0$;
\item $\eta_1=0$ and $\beta+ \eta_2\ge 0$.
\end{itemize}
\end{theo}
(In the case of $n=1$, then the Laplace--Beltrami operator is not defined, but obviously any closed linear operator $CD^{A,L}_0+D$  on $\partial X=\mathbb C^2$ is just a $2\times 2$-matrix and hence by Theorem~\ref{gener} no assumption has to be made on $\beta,\eta_1,\eta_2$ in order for $\mathbf A$ to be a generator. This turns out to be a special case of Theorem~\ref{gener-network} below.)

By construction the first coordinate of $(e^{t\mathbf A}{\mathbf u}_0)_{t\ge 0}$ solves the heat equation
\[
\frac{\partial u}{\partial t}(t,x)=\Delta u(t,x),\qquad t\ge 0,\ x\in \Omega,
\]
with dynamic boundary conditions $\rm(GWBC)$ (cf.\ Example~\ref{exa-01}) for all initial data
\[
{\mathbf u}_0:=\begin{pmatrix}u(0,\cdot)\\ f(0,\cdot)\end{pmatrix}\in L^2(\Omega)\times L^2(\partial \Omega).
\]


If $\beta>0$, then it is possible to associate a quadratic form with this operator -- boundedness (resp., uniform exponential stability) of $(e^{t\mathbf A})_{t\ge 0}$ can then be deduced from the accretivity of said form and the Lumer--Phillips theorem if $\eta_3=0$ (resp., if $\eta_3>0$). If $\beta=0$, this method is not applicable any more. The following holds by the argument presented in Remark~\ref{exa:misc}, though.

\begin{cor}
Under the assumptions of Theorem~\ref{theo:main-ww}, let $\beta= 0$ and additionally $\eta_1\ge 0$ or $\eta_2\ge 0$. If  $ \eta_3= 0$, then $(e^{t\mathbf A})_{t\ge 0}$ is bounded and in fact AAP, but not  uniformly exponentially stable; if $\eta_3< 0$, then  $(e^{t\mathbf A})_{t\ge 0}$ is uniformly exponentially stable.
\end{cor}

\begin{proof}
It suffices to observe that both $\Delta_{\partial \Omega}$ and $\DN_0$ are non-invertible negative semi-definite operators, hence they generate analytic and bounded semigroups (even AAP, since they have compact resolvent) that are not uniformly exponentially stable, and the same holds for their linear combination $D$ unless the bounded perturbation $\eta_3 \Id$ shifts its spectrum. Now, $A_0$ is the Laplacian with Dirichlet boundary conditions, hence a negative definite operator: Thus,~\eqref{eq:nonres-bis} is trivally satisfied.
\end{proof}
	
 It remains an open question whether boundedness or uniform exponential stability hold for $ \beta<0$, too, provided $\eta_1>0$. If $\eta_1=0$, then by~\cite[Thm.~9]{BanBelRei06} $\mathbf A$ does have eigenvalues in the open right halfplane.

\begin{rem}\label{rem:spectr}
Also some basic spectral properties of operators matrices can be investigated by the factorization methods described in Section~\ref{sec:nondiago}: The following observation goes back to~\cite{Eng99,KraMugNag03b}. 

Let $\lambda$ lie in the resolvent set of $A_0$. Then both the resolvent operator $R(\lambda,A_0)$ and the operator $D_\lambda^{A,L}$ introduced in Lemma~\ref{dirichdef} are well-defined and the factorization in~\eqref{klen1} can be refined as follows if $B=0$, as was observed in~\cite[\S~II.3]{Eng98}:
\begin{equation}\label{fac}
\begin{array}{rl}
\lambda-{\mathbf{A}}
\;\;=&\mathbf{L}_\lambda {\mathbf{D}}_\lambda
\mathbf{M}_\lambda\\
:=&\begin{pmatrix}
\Id & 0\\
-C R(\lambda,A_0) & \Id
\end{pmatrix}
\begin{pmatrix}
\lambda- A_0 & 0\\
0 & \lambda-CD_\lambda^{A,L}-D
\end{pmatrix}
\begin{pmatrix}
\Id & -D_\lambda^{A,L}\\
0 & \Id
\end{pmatrix}.
\end{array}
\end{equation}
Recall that by assumption $C$ is bounded from $[D(A_0)]$ to $\partial X$. Hence, the operators $\mathbf{L}_\lambda$, $\mathbf{M}_\lambda$ are isomorphism on $\mathbf X$, hence $\lambda-{\mathbf{A}}$ is invertible if and only if the diagonal matrix (with diagonal domain) ${\mathbf{D}}_\lambda$ is invertible. We conclude that $\lambda\in \sigma ({\mathbf{A}})$ if and only if $\lambda\in\sigma(D+CD^{A,L}_\lambda)$. 

The standing assumption that $\lambda$ is not a spectral value of $A_0$ is not satisfactory and can be proved to be actually unnecessary in the concrete application that is most relevant to us. 

Let us define for all $\mu\in \mathbb C$ the operator $\Delta_\mu$ as the Laplacian on $L^2(\Omega)$ with boundary conditions
\[
\frac{\partial u}{\partial \nu}+\mu u_{|_{\partial\Omega}}=0
\]
(observe that $\Delta_\mu$ is negative definite for $\mu\ge 0$). Then it is known that $\lambda$ is element of the spectrum of $\Delta_\mu$ if and only if $\mu$ is element of the spectrum of the Dirichlet-to-Neumann operator $\DN_\lambda$ introduced in Example~\ref{exa-01} (see e.g.~\cite[Thm.~3.1]{AreMaz12}); and in particular, $\lambda$ is element of the spectrum of $\Delta_\lambda$ if and only if $\lambda$ is element of the spectrum of $\DN_\lambda$. Indeed, this a classical result that goes back to Krein under the additional assumption that $\lambda$ lies in the resolvent set of the Laplacian with Dirichlet boundary conditions, but this condition was shown to be superfluous in~\cite{AreMaz12}, provided one consider a more general definition of $\DN_\lambda$: Let for \emph{any} $\lambda\in \mathbb C$
\[
C(\lambda):=\left\{(g,h)\in L^2(\partial \Omega)\times L^2(\partial \Omega):\exists u\in H^1(\Omega)\hbox{ s.t. }\Delta u=\lambda u,\ u_{|_{\partial\Omega}}=g,\ \frac{\partial u}{\partial \nu}=h \right\},
\]
\[
K(\lambda):=\left\{h\in L^2(\partial \Omega):\exists u\in H^1(\Omega)\hbox{ s.t. }\Delta u=\lambda u,\ u_{|_{\partial\Omega}}=0,\ \frac{\partial u}{\partial \nu}=h \right\},
\]
and
\[
L^2_\lambda(\partial \Omega):=L^2(\partial \Omega)\ominus K(\lambda).
\]
Then the (generalized) Dirichlet-to-Neumann operator $\DN_\lambda$ is defined as the operator in $L^2_\lambda(\partial \Omega)$ whose graph is $C(\lambda)\cap (L^2_\lambda(\partial \Omega)\times L^2_\lambda(\partial \Omega))$. (If $\lambda$ lies in the resolvent set of the Laplacian with Dirichlet boundary conditions, then $K(\lambda)$ is a singleton and therefore $L^2_\lambda(\partial\Omega)=L^2(\partial\Omega)$.)

Now, observe that by definition $\lambda$ is an eigenvalue of $\Delta_\lambda$ precisely if there exists a nontrivial solution of $\eqref{EBP}$, i.e., if and only if $\lambda $ is an eigenvalue of $\mathbf A$ with $\beta=1$ and $\eta_1=\eta_2=\eta_3=0$. In this case, one deduces again from~\cite[Thm.~3.1]{AreMaz12} that the dimension of the eigenspace of $\mathbf A$ associated with $\lambda$ agrees with the dimension of the null space of $\lambda-\DN_\lambda$.
\end{rem}

\subsection{A coupled system}
In the literature, dynamic boundary conditions are commonly imposed on the \emph{trace} of the solution of partial differential equations. However, one may also consider a coupled system
\begin{equation}\label{cenn1}
\left\{
\begin{array}{rcll}
\frac{\partial u}{\partial t}(t,x)&=&\Delta u(t,x)-p(x)u(t,x), &t\geq 0,\; x\in\Omega,\\
\frac{\partial w}{\partial t}(t,z)&=& \Delta_{\partial \Omega} w(t,z)-q(z)w(t,z), &t\geq 0,\;z\in\partial\Omega,\\
 w(t,z)&=&\frac{\partial u}{\partial \nu}(t,z), &t\geq 0,\; z\in\partial \Omega,\\
 u(0,x)&=&f(x), &x\in\Omega,\\
 w(0,z)&=&h(z), &z\in\partial\Omega,
\end{array}
\right.
\end{equation}
of the form of those studied in~\cite[\S~3]{CasEngNag03}. Here $\Omega$ is a bounded open domain of ${\mathbb R}^n$, $n\geq 2$, with smooth boundary $\partial\Omega$, and $p\in L^\infty(\Omega)$, $q\in L^\infty(\partial\Omega)$.

Set
$$X:=L^2(\Omega),\quad \partial X:=L^2(\partial \Omega).$$
Define the operator
$$Au:=\Delta u-pu
$$
with domain
$$ D(A):=\left\{ u\in H^1(\Omega):\Delta u\in L^2(\Omega)\hbox{ and }\frac{\partial u}{\partial \nu}\in L^2(\partial \Omega)\right\}\subset H^\frac{3}{2}(\Omega),$$
and let additionally
$$Lu:=\frac{\partial u}{\partial \nu},\qquad u\in D(L):= D(A),$$
$$B=C=0,$$
$$Dw:=\Delta_{\partial \Omega} w-qw,\qquad w\in D(D):=H^2(\partial\Omega),$$
i.e., $D$ is (up to a bounded perturbation) the Laplace--Beltrami operator on $\partial \Omega$. In order to satisfy the Assumptions~\ref{basic-2} we  have to set $D(C):=D(A)$, $\partial Y:=\partial X$. Then, $A_0=A_{|\ker L}$ is (up to a bounded perturbation) the Laplacian with Neumann boundary conditions, and one sees that the Assumptions~\ref{basic} and~\ref{basic-2} are satisfied, hence Theorem~\ref{gener}.(2) applies and we conclude that~\eqref{cenn1} is governed by an analytic semigroup on $L^2(\Omega)\times L^2(\partial\Omega)$, as was already shown in~\cite[\S~3]{CasEngNag03}. 

The following holds applying the arguments in Remark~\ref{exa:misc}.

\begin{prop}
Let $0\leq p\in L^\infty(\Omega)$, $0\leq q\in L^\infty(\partial\Omega)$. 
If $p(x)>0$ for a.e.\ $x\in \Omega$ \emph{or} $q(z)>0$ for a.e.\ $z\in \partial \Omega$, then $(e^{t\mathbf A})_{t\ge 0}$ is bounded and in fact AAP. If both $p(x)>0$ for a.e.\ $x\in \Omega$ \emph{and} $q(z)>0$ for a.e.\ $z\in \partial \Omega$, then $(e^{t\mathbf A})_{t\ge 0}$ is uniformly exponentially stable.
\end{prop}

Observe that if $p\equiv 0$ and $q\equiv 0$ the semigroups $(e^{tA_0})_{t\ge 0},(e^{tD})_{t\ge 0}$ are analytic and bounded, but the non-resonance condition~\eqref{eq:nonres-0} is not satisfied, since neither $A_0$ nor $D$ are invertible. 

\subsection{The Laplacian on a network}\label{exas3}
Finally, let us consider a diffusion equation on a network-like structure: This is a standard problem in the theory of so-called \emph{quantum graphs}, see e.g.~\cite{Bel94,BerKuc13,Mug14}. Let more precisely $\mG=(\mV,\mE)$ be a (possibly infinite) graph whose $\mV\times \mE$ incidence matrix we denote by $\mathcal I=(\iota_{\mv\me})$: this allows us to introduce the notation
\[
\mE_\mv:=\{\me\in \mE:\iota_{\mv\me}\neq 0\},\qquad \mv\in \mV
\] 
for the set of edges incident in $\mv$. Attach an interval of length $\rho_\me$ to each edge $\me$: We will assume for the sake of simplicity that $\rho_\me=1$ for all $\me\in \mE$, but all our results prevail if there exist $r,R>0$ such that merely $r<\rho_\me<R$ for all $\me\in \mE$. In this way $\mG$ turns into a \emph{metric graph}, in the sense of~\cite[Chapt.~3]{Mug14}. We define on it a Laplacian by
\[
\Delta u:=\left(\frac{d^2 u_\me}{dx^2}\right)_{\me \in \mE},\qquad
\]
on which one typically imposes in the nodes continuity
\begin{equation}\tag{Cc}\index{$\rm(Cc)$}\index{Condition!continuity}
u_\me(\mv)=u_\mf(\mv)=:u(\mv),\qquad \hbox{for all }\me,\mf\in \mE_\mv,\; \mv\in \mV,
\end{equation} 
along with Kirchhoff-type conditions
\begin{equation}\tag{Kc}\index{$\rm(Kc)$}
\partial_\nu u(\mv):=\sum_{\me\in \mE}\frac{\partial u_\me}{\partial \nu} (\mv)=0,\qquad \hbox{for all }\mv\in \mV,
\end{equation}
where $\frac{\partial u_\me}{\partial \nu} (\mv)=-u'_\me (0)$ if $\mv$ is the initial endpoint of the edge $\me$ and $\frac{\partial u_\me}{\partial \nu} (\mv)=u'_\me (1)$ if $\mv$ is the terminal endpoint of the edge $\me$ -- or equivalently, using the incidence matrix $\mathcal I$,
$\frac{\partial u_\me}{\partial \nu} (\mv)=\iota_{\mv\me}u'_\me(\mv)$. Under  $\rm(Cc)$, the vector of nodal values 
\[
u_{|\mV}\in {\mathbb C}^\mV
\]
of a function $u$ is well-defined. We fall within the scope of the theory developed in Section~\ref{sec:nondiago} if we pick a subset $\mV_0\subset \mV$ and replace there the conditions $\rm(Kc)$ by dynamic ones of the form
\begin{equation}\tag{dKc}
\frac{\partial u}{\partial t}(t,\mv)=-\partial_\nu u (t,\mv),\qquad \hbox{for all }\mv\in \mV_0,\ t\ge 0,
\end{equation}
along with
\begin{equation}\tag{Kc'}
\partial_\nu u(\mv)=0,\qquad \hbox{for all }\mv\in \mV\setminus \mV_0.
\end{equation}
This kind of boundary condition has a long history and interesting applications: We refer to~\cite{MugRom07} for a list of historical references related to this problem.
Generalizing this setting on the line of $\rm(GWBC)$ we are thus led to study the initial-boundary value problem
\begin{equation*}
\left\{
\begin{array}{rcll}
\frac{\partial u_\me}{\partial t}(t,x)&=& u''_\me(t,x), &t\ge 0,\; x\in(0,1), \me\in \mE,\\
u_\me(t,\mv)&=&u_\mf(t,\mv)=:u(t,\mv), &t\ge 0,\; \me,\mf\in \mE_\mv,\\
\frac{\partial u}{\partial t}(t,\mv)&=& -\beta_\mv \partial_\nu u(t,\mv)+\eta_\mv\DN_0 u(t,\mv)+\gamma_\mv u(t,\mv), &t\ge 0,\; \mv\in \mV_0,\\
 \partial_\nu u(t,\mv)&=&0, &t\ge 0,\; \mv\in \mV\setminus \mV_0,\\
u_\me(0,x)&=&u^0_\me(x),& x\in (0,1),\ \me\in \mE,\\
u_\mv(0)&=&f(\mv),& \mv\in \mV_0,
\end{array}
\right.
\end{equation*}
where $\beta,\eta,\gamma\in \mathbb R^\mV$ are general vectors on whose sign no assumption is made, and $\DN_0$ is the Dirichlet-to-Neumann map of the graph $\mG$ with respect to the node set $\mV_0$, i.e., for any $f,g\in \mathbb C^\mV$ and all $\lambda \in \mathbb C$ one defines $\DN_0f:=g$ if there exists $u\in H^2\big((0,1);\ell^2(\mE)\big)$ such that
\begin{equation}\label{bvp-dton-1dim}
\left\{
\begin{array}{rcll}
 u''_\me(x)&=&0,\qquad &x\in [0,1],\ \me\in \mE,\\
\partial_\nu u(\mv)&=&0, &\mv \in \mV\setminus\mV_0,\\
 \exists u(\mv)&:=&u_\me(\mv)=u_\mf(\mv),\qquad &\hbox{for all }\me,\mf\in \mE_\mv,\; \mv\in \mV,\\
-\partial_\nu u(\mv)&=&g(\mv), &\mv \in \mV_0.
\end{array}
\right.
\end{equation}
It is known that if $\mV_0=\mV$, then $\DN_0$ agrees with with the matrix $\mathcal I\mathcal I^T$, i.e., with the so-called discrete Laplacian of $\mG$, cf.~\cite[Chapter 2]{Mug14}.

Let us consider the Hilbert vector-valued space 
\[
X:= L^2\big((0,1);\ell^2(\mE)\big)
\]
as well as the weighted sequence spaces
\[\partial Y:=\{f\in \ell^2_{\rm \deg} (\mV_0):\mathcal I^T f\in \ell^2(\mE_0)\},\quad \partial X:=\ell^2_{\rm deg}(\mV_0),
\]
where the weight function $\deg$ is defined by
\[
\deg(\mv):=\sum_{\me\in \mE}|\iota_{\mv\me}|,\qquad \mv\in \mV_0,
\] -- equivalently, $\deg(\mv)$ is the cardinality of $\mE_\mv$ -- and $\mE_0$ is the set of all edges that are incident to at least one element of $\mV_0$.
Furthermore, we let $B=0$ and define the operators $A$, $C$, $D$, and $L$ by
$$Au:=\Delta u,\quad Cu:=-\beta\partial_\nu u_{|\mV_0},\quad Lu:=u_{|\mV_0},$$
with common domain
$$D(A)=D(C)=D(L):=\{u\in H^2\left((0,1);\ell^2(\mE):u\hbox{ satisfies $\rm(Cc)$ and $\rm(Kc')$}\right\},$$
along with
\[
Df:=\eta\DN_0 f+\gamma f,\qquad D(D):=\left\{
\begin{array}{ll}
\{f\in \partial Y:\DN_0 f \in \partial X\}\qquad &\hbox{if }\eta\neq 0,\\
\partial X &\hbox{otherwise}.
\end{array}
\right.
\]

With these choices $A_0$ becomes simply a family of second derivatives $X$ with Dirichlet boundary condition on all nodes that belong to $\mV_0$, while continuity and Kirchhoff-type conditioins are still imposed in the nodes that belong to $\mV\setminus\mV_0$. It is well-known that $A_0$ generates an analytic semigroup on $X$, see e.g.~\cite{Nic87} for the case of finite graphs and~\cite[Chapter~6]{Mug14} for the general case.

\begin{theo}\label{gener-network} Let $\gamma\in \ell^\infty(\mV)$. Then the following assertions hold.
\begin{enumerate}[(1)]
\item If $\deg\in \ell^\infty(\mV_0)$, then $\mathbf A$ generates an analytic semigroup on $X\times \partial X$.

\item If $\deg\not\in \ell^\infty(\mV_0)$, then $\mathbf A$ generates an analytic semigroup on $X\times \partial X$ if and only if 
\[
\mV_0^-:=\{\mv \in \mV_0:\beta_\mv+ \eta_\mv<0\}
\]
is a finite set.
\end{enumerate}
\end{theo}

In graph theoretical language, the condition $\deg\in \ell^\infty(\mV_0)$ is expressed by saying that the subgraph of $\mG$ induced by $\mV_0$ is uniformly locally finite.

In the case of finite $\mV_0$ (clearly, a special case of (1)), well-posedness has already been observed in~\cite[\S~19]{Bel94} and also in~\cite{MugRom07}.

\begin{proof}
The operator $A_0$ is certainly negative semidefinite as a simple integration by parts shows: It turns out that invertibility of $A_0$ is related to certain geometric properties of the graph. In order to avoid unnecessary technicalities, let us instead take some $\lambda>0$, which certainly lies in the resolvent set of $A_0$.

First of all, observe that 
\[
D^{A,L}_\lambda(\partial X)\subset H^2\big((0,1);\ell^2(\mE)\big)\hookrightarrow \left[[D(A_0)],X\right]_\alpha\qquad\hbox{ for all }\alpha<\frac{1}{4}.
\]
In view of Theorem~\ref{gener}.(1), $\mathbf A$ generates an analytic semigroup if and only if so do $A_0$ and $D+CD_\lambda^{A,L}$. Now, we already know that $A_0$ is a negative semi-definite operator, hence we can restrict our attention to the operator $CD_\lambda^{A,L}+D=(\beta+\eta)\DN_\lambda+\gamma \Id$.

Now, $\DN_\lambda$ is bounded on $\partial X$ if (and only if) $\deg\in \ell^\infty(\mV_0)$ (cf.~\cite[Chapter~6]{Mug14} or~\cite[Thm.~9.3]{HaeKelLen12}), and in this case $CD_\lambda^{A,L}+D$ generates an analytic semigroup and the assertion follows from Theorem~\ref{gener}.

If $\deg\not\in \ell^\infty(\mV_0)$ but $\mV_0^-$ is a finite set, then $CD_\lambda^{A,L}+D$ is only a compact perturbation of the analytic semigroup generator $\DN_\lambda$, hence again the operator matrix $\mathbf A$ associated with the problem generates an analytic semigroup. If however $\mV_0^-$ is infinite, and hence $(\beta+\eta)\DN_\lambda$ is not the compact perturbation of an analytic semigroup generator, then it is possible to find an infinite sequence of positive eigenvalues -- i.e., $(\beta+\eta)\DN_\lambda	$ is not an analytic semigroup generator, nor is $\mathbf A$ by Theorem~\ref{gener}.
\end{proof}

\bibliographystyle{alpha}
\bibliography{../../referenzen/literatur}
\end{document}